\documentclass[12pt]{article}
\usepackage{amsmath,amsfonts,amssymb,amsthm}
\input amssym.def
\begin{document}
\def \Z{\Bbb Z}
\def \C{\Bbb C}
\def \R{\Bbb R}
\def \Q{\Bbb Q}
\def \N{\Bbb N}
\def \bR{\bf R}
\def \D{{\cal{D}}}
\def \E{{\cal{E}}}
\def \Lie{{\cal{L}}}
\def \S{{\cal{S}}}
\def \Y{{\cal{Y}}}
\def \wt{{\rm wt}}
\def \tr{{\rm tr}}
\def \span{{\rm span}}
\def \Ga{{\rm \Gamma}}
\def \GL{{\rm GL}}
\def \gl{{\rm gl}}
\def \Res{{\rm Res}}
\def \End{{\rm End}\;}
\def \Ind {{\rm Ind}}
\def \Irr {{\rm Irr}}
\def \Aut{{\rm Aut}}
\def \Hom{{\rm Hom}}
\def \mod{{\rm mod}}
\def \ann{{\rm Ann}}
\def \ad{{\rm ad}}
\def \rank{{\rm rank}\;}

\def \<{\langle} 
\def \>{\rangle} 
\def \e{{\bf e}}
\def \f{{\bf f}}
\def \l{\lambda }
\def \g{{\frak{g}}}
\def \h{{\hbar}}
\def \k{{\frak{k}}}
\def \sl{{\frak{sl}}}

\def \be{\begin{equation}\label}
\def \ee{\end{equation}}
\def \bex{\begin{example}\label}
\def \eex{\end{example}}
\def \bl{\begin{lem}\label}
\def \el{\end{lem}}
\def \bt{\begin{thm}\label}
\def \et{\end{thm}}
\def \bp{\begin{prop}\label}
\def \ep{\end{prop}}
\def \br{\begin{rem}\label}
\def \er{\end{rem}}
\def \bc{\begin{coro}\label}
\def \ec{\end{coro}}
\def \bd{\begin{de}\label}
\def \ed{\end{de}}
\
\newcommand{\n}{\:^{\times}_{\times}\:}
\newcommand{\nno}{\nonumber}
\newcommand{\nord}{\mbox{\scriptsize ${\circ\atop\circ}$}}
\newtheorem{thm}{Theorem}[section]
\newtheorem{prop}[thm]{Proposition}
\newtheorem{coro}[thm]{Corollary}
\newtheorem{conj}[thm]{Conjecture}
\newtheorem{example}[thm]{Example}
\newtheorem{lem}[thm]{Lemma}
\newtheorem{rem}[thm]{Remark}
\newtheorem{de}[thm]{Definition}
\newtheorem{hy}[thm]{Hypothesis}
\makeatletter
\@addtoreset{equation}{section}
\def\theequation{\thesection.\arabic{equation}}
\makeatother
\makeatletter

\begin{center}{\Large \bf On certain generalizations of 
twisted affine Lie algebras
 and quasimodules for $\Gamma$-vertex algebras}
\end{center}

\begin{center}
{Haisheng Li\footnote{Partially supported by NSA grant
H98230-05-1-0018}\\
Department of Mathematical Sciences, Rutgers University, Camden, NJ 08102\\
and\\
Department of Mathematics, Harbin Normal University, Harbin, China}
\end{center}

\begin{abstract}
We continue a previous study on
$\Gamma$-vertex algebras and their quasimodules.
In this paper we refine certain known results and 
we prove that for any $\Z$-graded vertex algebra $V$ and a positive
integer $N$, the category of $V$-modules is naturally isomorphic to
the category of quasimodules of a certain type for $V^{\otimes N}$.
We also study certain generalizations of twisted affine Lie
algebras and we relate such Lie algebras to vertex algebras and their
quasimodules, in a way similar to that twisted affine Lie algebras are
related to vertex algebras and their twisted modules.
\end{abstract}

\section{Introduction}
It has been fairly well known that infinite-dimensional Lie algebras
such as (untwisted) affine Kac-Moody Lie algebras and the Virasoro Lie
algebra through their highest weight modules can be associated with
vertex operator algebras and modules (cf. \cite{fz}, \cite{dl}).  
On the other hand, it was known (see
\cite{flm}, \cite{li-twisted}) that twisted affine Lie algebras
through their highest weight modules can also be associated with
vertex operator algebras and their twisted modules.

In \cite{gkk}, Golenishcheva-Kutuzova and Kac introduced and studied a
notion of $\Gamma$-conformal algebra with $\Ga$ a group.  As it was
proved therein, a $\Ga$-conformal algebra structure on a vector space
$\g$ exactly amounts to a Lie algebra structure on $\g$ together with
a group action of $\Ga$ on $\g$ by automorphisms such that for any
$a,b\in \g$, $[ga,b]=0$ for all but finitely many $g\in \Ga$. To each
$\Ga$-conformal algebra $\g$, they associated an infinite-dimensional
Lie algebra whose underlying vector space is a certain quotient space
of $\g\otimes \C[t,t^{-1}]$.  Several families of infinite-dimensional
Lie algebras, including centerless twisted affine Lie algebras and
quantum torus Lie algebras (see [GKL1-2]), were realized in terms of
$\Ga$-conformal algebras, and some new infinite-dimensional Lie algebras
were also constructed.

In \cite{li-ga}, to associate vertex algebra-like structures to Lie algebras 
like quantum torus Lie algebras, we studied ``quasilocal'' subsets of
$\Hom (W,W((x)))$ for any given vector space $W$ and we proved that
any quasilocal subset generates a vertex algebra in a certain
canonical way. However, the vector space $W$ 
under the obvious action is {\em not} a module for the vertex algebras
generated by quasilocal subsets. Then a new notion of what we called a
quasimodule naturally arose.
For a vertex algebra $V$, a quasimodule is a vector space $W$
equipped with a linear map $Y_{W}$ from $V$ to $\Hom (W,W((x)))$
satisfying the condition that $Y_{W}({\bf 1},x)=1$ and for $u,v\in V$,
there exists a nonzero polynomial $p(x_{1},x_{2})$ such that 
\begin{eqnarray*}
& &x_{0}^{-1}\delta\left(\frac{x_{1}-x_{2}}{x_{0}}\right)
p(x_{1},x_{2})Y_{W}(u,x_{1})Y_{W}(v,x_{2})\\
& &\hspace{1cm}-x_{0}^{-1}\delta\left(\frac{x_{2}-x_{1}}{-x_{0}}\right)
p(x_{1},x_{2})Y_{W}(v,x_{2})Y_{W}(u,x_{1})\\
&=&x_{2}^{-1}\delta\left(\frac{x_{1}-x_{0}}{x_{2}}\right)
p(x_{1},x_{2})Y_{W}(Y(u,x_{0})v,x_{2}).
\end{eqnarray*}
In terms of this notion, any vector space $W$ is a quasimodule 
for the vertex algebras generated by quasilocal subsets.
Taking $W$ to be a highest weight module for
the quantum torus Lie algebra we obtain a vertex algebra
with $W$ as a quasimodule.

For a vertex algebra $V$, the notion of quasimodule is intrinsically related
to the notion of twisted module (with respect to a finite order
automorphism).  As we mentioned before, highest weight modules for
twisted affine Lie algebras are naturally twisted modules for the
vertex operator algebras associated to untwisted affine Lie algebras.
It was showed in \cite{li-ga} that highest weight modules of fixed
level for a twisted affine Lie algebra, which is viewed as an
invariant subalgebra of the untwisted affine Lie algebra, are
naturally quasimodules for the vertex algebras associated to an
untwisted affine Lie algebra.  Motivated by these two facts, in
\cite{li-tq}, by using a result of Barron, Dong and Mason \cite{bdm},
we established a canonical connection
between twisted modules and quasimodules 
for general vertex operator algebras.

Let $V$ be a vertex operator algebra and 
let $\sigma$ be an order-$N$ automorphism
of $V$. For a $\sigma$-twisted $V$-module $W$, the vertex operator map 
$Y_{W}$, a linear map from $V$ to $\Hom (W,W((x^{1/N})))$,
satisfies the following invariance property
\begin{eqnarray*}
Y_{W}(\sigma v,x)=\lim_{x^{1/N}\rightarrow \omega_{N}x^{1/N}}Y_{W}(v,x)
\;\;\;\mbox{ for }v\in V.
\end{eqnarray*}
Set $\Ga=\<\sigma\>\subset \Aut V$ and let $\phi: \Ga\rightarrow
\C^{\times}$
is the group embedding determined by $\phi(\sigma)=\exp(2\pi i/N)$.
We call a quasimodule $(W,Y_{W})$ for $V$ a {\em
  $(\Ga,\phi)$-quasimodule}
if for any $u,v\in V$, there exists a nonnegative integer $k$ such
that the quasi-Jacobi identity holds with 
$p(x_{1},x_{2})=(x_{1}^{N}-x_{2}^{N})^{k}$ and such that 
the following invariance property holds:
$$Y_{W}(gv,x)=Y_{W}(\phi(g)^{L(0)}v,\phi(g)x)
\;\;\;\mbox{ for }v\in V,\; g\in \Ga.$$
What was proved in \cite{li-tq} is that 
the category of (weak) $\sigma$-twisted $V$-modules
is naturally isomorphic to the category of $(\Ga,\phi)$-quasimodules for $V$.

In \cite{li-ga}, partially motivated by the notion of $\Ga$-conformal algebra
in \cite{gkk} we formulated and studied a notion of $\Ga$-vertex algebra.
For any group $\Ga$, a $\Ga$-vertex algebra 
can be equivalently defined as a vertex algebra $V$ 
equipped with two group homomorphisms
$$R: \Ga\rightarrow GL(V),\;\;\; \phi: \Ga\rightarrow \C^{\times}$$
such that $R_{g}{\bf 1}={\bf 1}$ for $g\in \Ga$ and
$$R_{g}Y(v,x)R_{g}^{-1}=Y(R_{g}v,\phi(g)^{-1}x)
\;\;\;\mbox{ for }g\in \Ga,\; v\in V.$$
Let $V$ be a $\Z$-graded vertex algebra with $L(0)$ the grading
operator and let $\Ga$ be a group of automorphisms of $\Z$-graded
vertex algebra $V$. Let $\phi: \Ga\rightarrow \C^{\times}$ be any group
homomorphism. Define $R_{g}=\phi(g)^{-L(0)}g$ for $g\in \Ga$.
Then $V$ becomes a $\Ga$-vertex algebra.

In this paper, we formulate and study a notion of
quasimodule for $\Ga$-vertex algebras.
For a $\Ga$-vertex algebra $V$,
a $V$-quasimodule $W$ is a quasimodule for $V$ viewed as a vertex
algebra such that for any $u,v\in V$, the quasi-Jacobi identity holds
with $p(x_{1},x_{2})=(x_{1}-\alpha_{1}x_{2})\cdots
(x_{1}-\alpha_{r}x_{2})$ for some $\alpha_{1},\dots,\alpha_{r}\in
\phi(\Ga)\subset \C^{\times}$ and such that
\begin{eqnarray*}
Y_{W}(R_{g}v,x)=Y_{W}(v,\phi(g)x)
\;\;\;\mbox{ for }v\in V,\; g\in \Ga.
\end{eqnarray*}

Note that for any vector space $W$, the group $\C^{\times}$
acts on the space $\Hom(W,W((x)))$ by
$$R_{\lambda}a(x)=a(\lambda x)\;\;\;\mbox{ for }
a(x)\in \Hom (W,W((x))),\;\lambda\in \C^{\times}.$$
Let $\Ga$ be a subgroup of $\C^{\times}$.
A subset $S$ of $\Hom (W,W((x)))$ is said to be $\Ga$-local (see
\cite{gkk}) if for any $a(x),b(x)\in S$, there there exists
$\alpha_{1},\dots,\alpha_{r}\in \Ga$ such that
$$(x_{1}-\alpha_{1}x_{2})\cdots
(x_{1}-\alpha_{r}x_{2})[a(x_{1}),b(x_{2})]=0.$$
By refining a result of \cite{li-ga}, we prove that
every $\Ga$-local subset of $\Hom (W,W((x)))$ generates 
a $\Ga$-vertex algebra with $W$ as a quasimodule.
We also obtain an analogue of Borcherds' commutator formula
for quasimodules.

A conceptual result of Barron, Dong and Mason \cite{bdm} is that for
any vertex operator algebra $V$ and for any positive integer $N$, the
category of (weak) $V$-modules is canonically isomorphic to the
category of (weak) $\sigma$-twisted $V^{\otimes N}$-modules, where
$\sigma$ is a permutation automorphism of $V^{\otimes N}$.  In this
paper we show that for any $\Z$-graded vertex algebra $V$ and for any
positive integer $N$, a $V$-module structure on a vector space $W$
exactly amounts to a quasimodule structure for $V^{\otimes N}$ viewed
as a $\Ga$-vertex algebra with $\Ga=\<(12\cdots N)\>$, where
$(12\cdots N)$ denotes the permutation automorphism of $V^{\otimes
N}$.  This result can be considered as a version of
Barron, Dong and Mason's theorem in terms of quasimodules.

As we mentioned earlier, 
Golenishcheva-Kutuzova and Kac gave a construction of
infinite-dimensional Lie algebras from $\Ga$-conformal algebras
(in the sense of \cite{gkk}).
In this paper, we slightly extend 
their construction with a central extension being included.
Let $\g$ be any Lie algebra equipped with
a symmetric invariant bilinear form $\<\cdot,\cdot\>$.
Associated to the pair $(\g,\<\cdot,\cdot\>)$, one has the
untwisted affine Lie algebra 
$\hat{\g}=\g\otimes \C[t,t^{-1}]\oplus \C {\bf k}$
(see \cite{kac1}). 
Let $\Ga$ be a group of
automorphisms of $\g$, preserving the bilinear form, such that
for any $a,b\in \g$,
$$[ga,b]=0,\ \ \ \ \<ga,b\>=0\;\;\;
\mbox{ for all but finitely many }g\in \Ga,$$
let $\phi: \Ga\rightarrow \C^{\times}$ be any group homomorphism.
We construct a Lie algebra $\hat{\g}[\Ga]$ as a quotient space
of $\hat{\g}$. Furthermore, we prove that the category of
``restricted'' $\hat{\g}[\Ga]$-modules of
level $\ell$ is canonically isomorphic to
the category of quasimodules for $V_{\hat{\bar{\g}}}(\ell,0)$ viewed as a
$\Ga$-vertex algebra, where $\bar{\g}$ is a certain Lie algebra with
$\bar{\g}=\g$ for $\phi$ an injective homomorphism.
We also extend Golenishcheva-Kutuzova and Kac's notion of
$\Ga$-conformal algebra to include higher order singularity.

Note that extended affine Lie algebras (cf. [S1,2], \cite{my}, \cite{aabgp})
form a relatively large family of Lie algebras, 
including finite-dimensional simple Lie
algebras, (twisted and twisted) affine Lie algebras,
toroidal Lie algebras, and quantum torus algebras.
As these special extended affine Lie algebras have been associated with
vertex algebras (\cite{fz}, \cite{bbs}, \cite{li-ga}),
our naive hope is that every extended affine Lie algebra
can be realized as a generalized twisted affine Lie algebra $\hat{\g}[\Ga]$, 
so that all the extended affine Lie
algebras can be associated with vertex algebras and quasimodules.

This paper is organized as follows: In Section 2, we reformulate
and refine certain results on $\Gamma$-vertex algebras and quasimodules.
In Section 3, we give a natural isomorphism between the category of
$V$-modules and a certain subcategory of quasimodules for $V^{\otimes N}$.
In Section 4, we study certain generalizations of 
twisted affine Lie algebras and we relate them with vertex algebras 
and quasimodules.

\section{$\Gamma$-vertex algebras and their quasimodules}
In this section we reformulate the notion of
$\Gamma$-vertex algebra of \cite{li-ga} and we define a notion of quasimodule
for a $\Ga$-vertex algebra. Certain results of
\cite{li-ga} are refined and an analogue of
Borcherds' commutator formula is obtained.

We first recall from \cite{li-ga} the notion of quasimodule 
for a vertex algebra $V$, which generalizes
the notion of module. A {\em $V$-quasimodule}
is a vector space $W$ equipped with a linear map 
$$Y_{W}: V\rightarrow \Hom (W,W((x)))\subset (\End W)[[x,x^{-1}]]$$
satisfying the condition that $Y_{W}({\bf 1},x)=1$ and 
for $u,v\in V$, 
there exists a nonzero polynomial $p(x_{1},x_{2})$ such that
\begin{eqnarray}\label{epjacobi}
& &x_{0}^{-1}\delta\left(\frac{x_{1}-x_{2}}{x_{0}}\right)
p(x_{1},x_{2})Y_{W}(u,x_{1})Y_{W}(v,x_{2})\nonumber\\
& &\hspace{2cm}-x_{0}^{-1}\delta\left(\frac{x_{2}-x_{1}}{-x_{0}}\right)
p(x_{1},x_{2})Y_{W}(v,x_{2})Y_{W}(u,x_{1})\nonumber\\
& &=x_{2}^{-1}\delta\left(\frac{x_{1}-x_{0}}{x_{2}}\right)
p(x_{1},x_{2})Y_{W}(Y(u,x_{0})v,x_{2}).
\end{eqnarray}

\bl{lsimple2}
Let $V$ be a vertex algebra. A $V$-quasimodule $(W,Y_{W})$ is
a $V$-module if and only if for $u,v\in V$, $Y_{W}(u,x)$ and
$Y_{W}(v,x)$ are mutually local in the sense that
$$(x_{1}-x_{2})^{k}Y_{W}(u,x_{1})Y_{W}(v,x_{2})=
(x_{1}-x_{2})^{k}Y_{W}(v,x_{2})Y_{W}(u,x_{1})$$
for some nonnegative integer $k$, depending on $u$ and $v$.
\el

\begin{proof} We only need to prove the ``if'' part.
For $u,v\in V$, there exists a nonzero polynomial
$p(x_{1},x_{2})$ such that (\ref{epjacobi}) holds. 
Using the delta-function substitution we have
\begin{eqnarray}\label{emult-jac}
& &p(x_{0}+x_{2},x_{2})\nonumber\\
& &\cdot\left(x_{0}^{-1}\delta\left(\frac{x_{1}-x_{2}}{x_{0}}\right)
Y_{W}(u,x_{1})Y_{W}(v,x_{2})
-x_{0}^{-1}\delta\left(\frac{x_{2}-x_{1}}{-x_{0}}\right)
Y_{W}(v,x_{2})Y_{W}(u,x_{1})\right)\nonumber\\
&=&p(x_{0}+x_{2},x_{2})x_{2}^{-1}\delta\left(\frac{x_{1}-x_{0}}{x_{2}}\right)
Y_{W}(Y(u,x_{0})v,x_{2}).
\end{eqnarray}
As $Y_{W}(u,x)$ and
$Y_{W}(v,x)$ are mutually local, 
the second factor of the left-hand side of
(\ref{emult-jac}) involves only finitely many 
negative powers of $x_{0}$. 
Recall from \cite{li-ga} the natural field-embedding
$$\iota_{x_{1},x_{2}}: \C(x_{1},x_{2})\rightarrow
\C((x_{1}))((x_{2})),$$
where $\C(x_{1},x_{2})$ denotes the field of rational functions.
Then we multiply both sides of  (\ref{emult-jac}) by
$\iota_{x_{2},x_{0}}(1/p(x_{0}+x_{2},x_{2}))$, obtaining
the usual Jacobi identity. Thus $(W,Y_{W})$ is a $V$-module.
\end{proof}

Let $\Gamma$ be a group which is fixed throughout this section, and
denote by $\C^{\times}$ the group of nonzero complex numbers.

\bd{dGphistructure}
{\em A {\em $\Gamma$-vertex algebra} is a vertex algebra $V$ equipped with 
group homomorphisms
\begin{eqnarray}
& &R: \Gamma\rightarrow \GL(V); \; g\mapsto R_{g}\\
& &\phi: \Gamma \rightarrow \C^{\times}
\end{eqnarray}
 such that $R_{g}({\bf 1})={\bf 1}$ for $g\in \Gamma$  and 
\begin{eqnarray}\label{e-abstract-conjugation-def}
R_{g}Y(v,x)R_{g}^{-1}=Y(R_{g}(v),\phi(g)^{-1}x)
\;\;\;\mbox{ for }g\in \Gamma,\; v\in V.
\end{eqnarray}}
\ed

In view of (\cite{li-ga}, Theorem 6.5), this notion is
equivalent to the notion of $\Ga$-vertex algebra defined in \cite{li-ga}.

\bex{example-Zva}
{\em Let $V$  be a $\Z$-graded vertex algebra
in the sense that $V$ is a vertex algebra equipped with a $\Z$-grading
$V=\coprod_{n\in \Z}V_{(n)}$ such that
$$u_{k}V_{(n)}\subset V_{(m+n-k-1)} \;\;\;
\mbox{ for }u\in V_{(m)},\; m,n,k\in \Z.$$
Denote by $L(0)$ the grading operator, i.e.,
$$L(0)v=nv\;\;\;\mbox{ for }v\in V_{(n)},\; n\in \Z.$$
Let $\Gamma$ be a group of grading-preserving 
automorphisms of $V$ and let $\phi$ be any group homomorphism {}from 
$\Gamma$ to $\C^{\times}$. Define $R: \Gamma\rightarrow \GL(V)$ by
\begin{eqnarray}
R_{g}(v)=\phi(g)^{-L(0)}(gv)\;\;\;\mbox{ for }g\in \Gamma,\; v\in V.
\end{eqnarray}
{}From \cite{li-ga}, $V$ becomes a $\Gamma$-vertex algebra.}
\eex

\br{rfaithful}
{\em Let $V$ be a $\Ga$-vertex algebra
and let $g\in \ker R$. We have $Y(v,x){\bf 1}=Y(v,\phi(g)^{-1}x){\bf 1}$ 
for $v\in V$. Thus $e^{x\D}v=e^{\phi(g)^{-1}x\D}v$.
If $\D\ne 0$, or equivalently, if
$V$ is not a classical commutative associative algebra,
we have $\phi(g)=1$. That is, if $\D\ne 0$,
we have $\ker R\subset \ker \phi$, so
we can replace $\Gamma$ by the quotient group $\Gamma/\ker R$ 
with a faithful action on $V$.}
\er

\bd{dGphi-quasimodule}
{\em Let $V$ be a $\Gamma$-vertex algebra. 
A {\em $V$-quasimodule} is
a quasimodule $(W,Y_{W})$ for $V$ viewed as a vertex algebra, 
satisfying the condition that 
\begin{eqnarray}
Y_{W}(R_{g}v,x)=Y_{W}(v,\phi(g)x)\;\;\;\mbox{ for }g\in \Gamma,\; v\in V,
\end{eqnarray}
and for $u,v\in V$, there exist 
$\alpha_{1},\dots,\alpha_{k}\in \phi(\Gamma)\subset \C^{\times}$
such that
\begin{eqnarray}
(x_{1}-\alpha_{1}x_{2})\cdots (x_{1}-\alpha_{k}x_{2})
[Y_{W}(u,x_{1}),Y_{W}(v,x_{2})]=0.
\end{eqnarray}}
\ed

\br{rtwisted-quasi}
{\em Let $V$ be a vertex operator algebra in the sense of \cite{flm}
and let $\sigma$ be an order-$N$ automorphism of $V$. Set $\Ga=\<\sigma\>$
and let $\phi$ be the group homomorphism from $\Ga$ to $\C^{\times}$,
defined by $\phi(\sigma)=\exp(2\pi i/N)$. Consider
$V$ as a $\Gamma$-vertex algebra as in Example \ref{example-Zva}.
It was proved in \cite{li-tq} that
the category of weak $\sigma$-twisted $V$-modules is isomorphic to
the category of quasimodules for $V$ viewed as a $\Ga$-vertex algebra. }
\er

Now, let $W$ be a general vector space. Set 
$$\E(W)=\Hom (W,W((x)))\subset (\End W)[[x,x^{-1}]].$$
A subset $S$ of $\E(W)$ is said to be {\em quasi-local}
if for any $a(x),b(x)\in S$, there exists a nonzero polynomial 
$p(x_{1},x_{2})$ such that
\begin{eqnarray}\label{equasi-locality}
p(x_{1},x_{2})a(x_{1})b(x_{2})=p(x_{1},x_{2})b(x_{2})a(x_{1}).
\end{eqnarray}
Let $\Gamma$ be a subgroup of $\C^{\times}$.
Following \cite{gkk}, we say $a(x),b(x)\in \E(W)$
are {\em $\Gamma$-local} if 
\begin{eqnarray}\label{egamma-locality}
(x_{1}-\alpha_{1} x_{2})\cdots (x_{1}-\alpha_{r}x_{2})a(x_{1})b(x_{2})
=(x_{1}-\alpha_{1} x_{2})\cdots (x_{1}-\alpha_{r}x_{2})b(x_{2})a(x_{1})
\end{eqnarray}
for some $\alpha_{1},\dots,\alpha_{r}\in \Gamma$.
A subset $S$ of  $\E(W)$
is said to be {\em $\Gamma$-local} if every pair in $S$ is $\Gamma$-local.
 
\br{rmodule-action}
{\em Note that $\E(W)$ is naturally a vector space over the field $\C((x))$. 
On the other hand, we define a group action of $\C^{\times}$ on $\E(W)$ 
with $\alpha\in \C^{\times}$ acting as $R_{\alpha}$ by
\begin{eqnarray}
R_{\alpha}a(x)=a(\alpha x)\;\;\;\mbox{ for }a(x)\in \E(W)
\end{eqnarray}
(cf. \cite{gkk}). Notice that for any quasi-local subset $S$ of $\E(W)$,
the subspace spanned by $R_{\alpha}(S)$ for $\alpha\in \C^{\times}$ is 
also quasi-local. Consequently, every maximal quasi-local subspace of 
$\E(W)$ is closed under the action of $\C^{\times}$.}
\er

Assume that $a(x),b(x)\in \E(W)$ are quasi-local. 
Notice that the commutativity relation (\ref{equasi-locality}) implies
$$p(x_{1},x_{2})a(x_{1})b(x_{2})\in \Hom (W,W((x_{1},x_{2}))).$$
We define
\begin{eqnarray*}
Y_{\E}(a(x),x_{0})b(x)
&=& \iota_{x,x_{0}}(1/p(x+x_{0},x))\Res_{x_{1}}
x_{1}^{-1}\delta\left(\frac{x+x_{0}}{x_{1}}\right)
\left(p(x_{1},x)a(x_{1})b(x)\right)\\
&=&\iota_{x,x_{0}}(1/p(x+x_{0},x))
\left(p(x_{1},x)a(x_{1})b(x)\right)|_{x_{1}=x+x_{0}}.
\end{eqnarray*}
Write
$$Y_{\E}(a(x),x_{0})b(x)=\sum_{n\in \Z}a(x)_{n}b(x) x_{0}^{-n-1}.$$
A quasilocal subspace $U$ of $\E(W)$ is said to be {\em closed} (under
$Y_{\E}$) if
\begin{eqnarray}
a(x)_{n}b(x)\;\;\;\;\mbox{ for }a(x),b(x)\in U,\; n\in \Z.
\end{eqnarray}

\br{ryalpha}
{\em For every $\alpha\in \C^{\times}$,
  $Y_{\E}^{(\alpha)}(a(x),x_{0})b(x)$
was defined in \cite{li-ga}, where
$$Y_{\E}^{(\alpha)}(a(x),x_{0})b(x)
=\iota_{x,x_{0}}(1/p(\alpha x+x_{0},x))
\left(p(x_{1},x)a(x_{1})b(x)\right)|_{x_{1}=\alpha x+x_{0}}.$$
It was proved (\cite{li-ga}, Proposition 3.9) that
$$Y_{\E}^{(\alpha)}(a(x),x_{0})b(x)
=Y_{\E}(a(\alpha x),\alpha^{-1}x_{0})b(x).$$
It is clear that a quasilocal subspace $U$ of $\E(W)$ is
closed under all the operations $Y_{\E}^{(\alpha)}$ for $\alpha\in
\Gamma$ (a subgroup of
$\C^{\times}$) if and only if $U$ is closed under $Y_{\E}$ and closed
under the action of $\Gamma$.}
\er

\bt{tga-old} 
Let $W$ be a vector space over $\C$
and let $\Gamma$ be a subgroup of $\C^{\times}$. 
For any $\Gamma$-local subset $S$ of 
$\E(W)$, there exist closed $\Gamma$-local subspaces $K$
of $\E(W)$ with the following property
\begin{eqnarray}
\{1_{W}\} \cup S\subset K,\;\;\;
R_{\alpha}a(x)=a(\alpha x)\in K
\;\;\;\mbox{ for }\alpha\in \Gamma,\; a(x)\in K,
\end{eqnarray}
among which the smallest subspace is denoted by 
$\<S\>_{\Gamma}$.
Furthermore, let $V$ be any closed $\Gamma$-local subspace of $\E(W)$ such that
$$1_{W}\in V\;\mbox{ and }\; R_{\alpha}a(x)=a(\alpha x)\in V
\;\;\;\mbox{ for }\alpha\in \Gamma,\; a(x)\in V.$$
Then $(V,1_{W},Y_{\E})$ carries the structure of a $\Gamma$-vertex algebra 
with $W$ as a quasimodule where $Y_{W}(a(x),x_{0})=a(x_{0})$ 
for $a(x)\in V$.
\et

\begin{proof} By  (\cite{li-ga}, Corollary 4.7), 
there exists a (unique) smallest $\Gamma$-local subspace
$\<S\>_{\Gamma}$, that contains $1_{W}$ and $S$ and that is closed
under $Y_{\E}^{(\alpha)}$ for $\alpha\in \Gamma$.
{}From Remark \ref{ryalpha}, $\<S\>_{\Gamma}$ is also the smallest
$\Gamma$-local subspace that is closed under $Y_{\E}$
and under the action of $\Gamma$.
This proves the first assertion.

For the second assertion, with $V$ being $\Gamma$-local, $V$ is quasilocal. 
By Theorems 6.3 and 6.5 of \cite{li-ga}, $(V,1_{W},Y_{\E})$  
carries the structure of a $\Gamma$-vertex algebra with $W$ as a quasimodule
where $Y_{W}(\alpha(x),x_{0})=\alpha(x_{0})$. 
For $\alpha\in \Gamma,\; a(x)\in V$, we have
$$Y_{W}(R_{\alpha}a(x),x_{0})=Y_{W}(a(\alpha x),x_{0})=a(\alpha x_{0})
=Y_{W}(a(x),\alpha x_{0}).$$
Thus, $(W,Y_{W})$ is a $V$-quasimodule.
\end{proof}

\br{rcomments}
{\em Notice that for $\Gamma=\{1\}$ (the trivial group), 
$\Ga$-locality becomes the usual locality (cf. \cite{li-local}). 
Assume that $a(x),b(x)\in \E(W)$ are local, i.e., 
there exists a nonnegative integer $k$ such that
$$(x_{1}-x_{2})^{k}a(x_{1})b(x_{2})=(x_{1}-x_{2})^{k}b(x_{2})a(x_{1}).$$
Then
\begin{eqnarray*}
& &Y_{\E}(a(x),x_{0})b(x)\nonumber\\
&=&x_{0}^{-k}\Res_{x_{1}}
x_{1}^{-1}\delta\left(\frac{x+x_{0}}{x_{1}}\right)
\left((x_{1}-x)^{k}a(x_{1})b(x)\right)\\
&=&\Res_{x_{1}}
x_{0}^{-1}\delta\left(\frac{x_{1}-x}{x_{0}}\right)
x_{0}^{-k}\left((x_{1}-x)^{k}a(x_{1})b(x)\right)\\
& &-\Res_{x_{1}}
x_{0}^{-1}\delta\left(\frac{x-x_{1}}{-x_{0}}\right)
x_{0}^{-k}\left((x_{1}-x)^{k}a(x_{1})b(x)\right)\\
&=&\Res_{x_{1}}
x_{0}^{-1}\delta\left(\frac{x_{1}-x}{x_{0}}\right)
(x_{1}-x)^{-k}\left((x_{1}-x)^{k}a(x_{1})b(x)\right)\\
& &-\Res_{x_{1}}
x_{0}^{-1}\delta\left(\frac{x-x_{1}}{-x_{0}}\right)
(-x+x_{1})^{-k}\left((x_{1}-x)^{k}b(x)a(x_{1})\right)\\
&=&\Res_{x_{1}}\left(
x_{0}^{-1}\delta\left(\frac{x_{1}-x}{x_{0}}\right)
a(x_{1})b(x)
-x_{0}^{-1}\delta\left(\frac{x-x_{1}}{-x_{0}}\right)
b(x)a(x_{1})\right),
\end{eqnarray*}
which agrees with the definition given in \cite{li-local}. 
Also, for $\Gamma=\{1\}$, a quasimodule is simply a module.
Thus Theorem \ref{tga-old} generalizes the corresponding result of 
\cite{li-local}. }
\er

\bp{phomo-rep}
Let $V$ be a $\Ga$-vertex algebra and let $W$ be a vector space
equipped with a linear map $Y_{W}$ from $V$ to $\E(W)$ such that
$Y_{W}({\bf 1},x)=1_{W}$.
Set
$$\overline{V}=\{ Y_{W}(v,x)\;|\; v\in V\}\subset \E(W).$$
Then $(W,Y_{W})$ carries the structure of a $V$-quasimodule if and only if
$\overline{V}$ is $\phi(\Ga)$-local, closed, and the map $Y_{W}$ 
is a homomorphism of vertex algebras 
{}from $V$ to $(\overline{V},1_{W},Y_{\E})$ such that
$$Y_{W}(R_{g}v,x)=Y_{W}(v,\phi(g)x)\;\;\;\mbox{ for }g\in \Ga.$$
\ep

\begin{proof} It was proved (\cite{li-ga}, Proposition 5.4) that
$(W,Y_{W})$ carries the structure of a quasimodule for $V$ viewed as a
  vertex algebra if and only if
$\overline{V}$ is quasi-local and closed,  and $Y_{W}$ is a
vertex-algebra homomorphism from $V$ to $(\overline{V},1_{W},Y_{\E})$.
Then it follows immediately.
\end{proof}

\bp{preduction}
Let $V$ be a $\Ga$-vertex algebra. Denote by $J$ the ideal of $V$ generated by
the elements $R_{h}v-v$ for $h\in \ker \phi\subset \Ga,\; v\in V$ and set $\overline{V}=V/J$ 
(the quotient vertex algebra). Then $\ker \phi$ acts trivially 
on $\overline{V}$
and $\overline{V}$ is a $\Ga/\ker \phi$-vertex algebra with
injective group homomorphism from $\Ga/\ker \phi$ to $\C^{\times}$.
Furthermore, for any $V$-quasimodule $(W,Y_{W})$,
$Y_{W}(u,x)=0$ for $u\in J$ and $W$ is naturally a $\overline{V}$-quasimodule.
\ep

\begin{proof} For $g\in \Ga,\; h\in \ker \phi,\; v\in V$, we have
$$R_{g}(R_{h}v-v)=R_{ghg^{-1}}R_{g}v-R_{g}v\in J,$$
as $ghg^{-1}\in \ker \phi$. {}From \cite{ll}, $J$ is linearly spanned by
all the coefficients of the formal series
$$Y(u,x)(R_{h}v-v),\;\;\; Y(R_{h}v-v,x)u$$
for $h\in \ker \phi,\; u,v\in V$.
As $R_{g}Y(a,x)b=Y(R_{g}a,\phi(g)x)R_{g}b$ for $g\in \Ga,\; a,b\in V$,
it follows that $J$ is stable under the action of $\Ga$.
{}From definition, $\ker \phi$ acts trivially on $\overline{V}$. 
Then the first assertion is clear. 
Now let $(W,Y_{W})$ be a $V$-quasimodule.
For $h\in \ker \phi,\; v\in V$, we have
$$Y_{W}(R_{h}v,x)=Y_{W}(v,\phi(h)x)=Y_{W}(v,x).$$
Thus $R_{h}v-v\in \ker Y_{W}$ for $h\in \ker \phi,\; v\in V$. 
As $Y_{W}$ is a homomorphism of vertex algebras (Proposition \ref{phomo-rep}),
$\ker Y_{W}$ is an ideal of $V$.
Consequently, $J\subset \ker Y_{W}$. Then
$W$ is naturally a $\overline{V}$-quasimodule.
\end{proof}

For quasimodules we have the following analogue of Borcherds'
commutator formula:

\bp{pmore1}
Let $V$ be a $\Ga$-vertex algebra and 
let $(W,Y_{W})$ be a $V$-quasimodule.
Let $\psi: \phi(\Ga)\rightarrow \Ga$ be any section of $\phi$. 
For $u,v\in V$, we have
\begin{eqnarray}
[Y_{W}(u,x_{1}),Y_{W}(v,x_{2})]=\sum_{\alpha\in \phi(\Ga)}\Res_{x_{0}}
x_{1}^{-1}\delta\left(\frac{\alpha x_{2}+x_{0}}{x_{1}}\right)
Y_{W}(Y(R_{\psi(\alpha)}u,\alpha^{-1}x_{0})v,x_{2}),
\end{eqnarray}
which is a finite sum. 
\ep

\begin{proof} Set $\overline{V}=\{Y_{W}(v,x)\;|\; v\in V \}\subset
  \E(W)$. By Proposition \ref{phomo-rep}, $\overline{V}$ 
is $\phi(\Ga)$-local and
closed and $Y_{W}$ is a vertex-algebra homomorphism
{}from $V$ to $(\overline{V},1_{W},Y_{\E})$.
Recall {}from \cite{li-ga} (Corollary 3.12) that
if $a(x),b(x)\in \E(W)$ satisfy the following relation
$$(x_{1}-\alpha_{1}x_{2})^{k_{1}}\cdots
(x_{1}-\alpha_{r}x_{2})^{k_{r}}
[a(x_{1}),b(x_{2})]=0,$$
where $\alpha_{1},\dots,\alpha_{r}$ are distinct nonzero complex
numbers and $k_{1},\dots,k_{r}$ are nonnegative integers,
then
$$Y_{\E}(a(\beta x),x_{0})b(x)\in \E(W)[[x_{0}]]\;\;\;
\mbox{ for }\beta\notin \{ \alpha_{1},\dots,\alpha_{r}\}$$
 and
\begin{eqnarray}
[a(x_{1}),b(x_{2})]=\sum_{\beta\in \C^{\times}}\Res_{x_{0}}
x_{1}^{-1}\delta\left(\frac{\beta x_{2}+x_{0}}{x_{1}}\right)
Y_{\E}(a(\beta x_{2}),\beta^{-1}x_{0})b(x_{2}).
\end{eqnarray}
(Note that $Y_{\E}^{(\beta)}(a(x),x_{0})b(x)=Y_{\E}(a(\beta x),\beta^{-1}x_{0})b(x)$.)
For $u,v\in V$, using the fact that $Y_{W}$ 
is a vertex-algebra homomorphism, we get
\begin{eqnarray*}
& &[Y_{W}(u,x_{1}),Y_{W}(v,x_{2})]\\
&=&\sum_{\beta\in \phi(\Ga)}\Res_{x_{0}}
x_{1}^{-1}\delta\left(\frac{\beta x_{2}+x_{0}}{x_{1}}\right)
\left(Y_{\E}(Y_{W}(u,\beta
x),\beta^{-1}x_{0})Y_{W}(v,x)\right)|_{x=x_{2}}\\
&=&\sum_{\beta\in \phi(\Ga)}\Res_{x_{0}}
x_{1}^{-1}\delta\left(\frac{\beta x_{2}+x_{0}}{x_{1}}\right)
\left(Y_{\E}(Y_{W}(R_{\psi(\beta)}u,
x),\beta^{-1}x_{0})Y_{W}(v,x)\right)|_{x=x_{2}}\\
&=&\sum_{\beta\in \phi(\Ga)}\Res_{x_{0}}
x_{1}^{-1}\delta\left(\frac{\beta x_{2}+x_{0}}{x_{1}}\right)
Y_{W}\left(Y(R_{\psi(\beta)}u,\beta^{-1}x_{0})v,x_{2}\right),
\end{eqnarray*}
which is a finite sum. 
\end{proof}

\bc{cfinite-issue}
Let $V$ be a $\Ga$-vertex algebra with a faithful action of $\Ga$.
Suppose that $V$ has a faithful quasimodule. 
Then $\Ga$ must be abelian and 
$\phi: \Ga\rightarrow \C^{\times}$ is injective. 
Furthermore, for any $u,v\in V$,
\begin{eqnarray}
Y(R_{g}u,x)v\in V[[x]]
\;\;\;\mbox{ for all but finitely many }g\in \Ga,
\end{eqnarray}
or equivalently,
\begin{eqnarray}
[Y(R_{g}u,x_{1}),Y(v,x_{2})]=0
\;\;\;\mbox{ for all but finitely many }g\in \Ga.
\end{eqnarray}
\ec

\begin{proof} Let $(W,Y_{W})$ be a faithful $V$-quasimodule.
Let $g\in \ker \phi$. Then
$$Y_{W}(R_{g}v,x)=Y_{W}(v,x)\;\;\;\mbox{ for all }v\in V.$$
It follows that $R_{g}=1$. As $\Ga$ acts faithfully on $V$, we have $g=1$.
This proves that $\phi$ is injective. Consequently,
$\Ga$ is abelian.

For $u,v\in V,\; g\in \Ga$, we have
\begin{eqnarray*}
Y_{W}(Y(R_{g}u,x_{0})v,x)=Y_{\E}(Y_{W}(R_{g}u,x),x_{0})Y_{W}(v,x)
=Y_{\E}(Y_{W}(u,\phi(g)x),x_{0})Y_{W}(v,x).
\end{eqnarray*}
As it was mentioned in the proof of Proposition \ref{pmore1},
$$Y_{\E}(Y_{W}(u,\phi(g)x),x_{0})Y_{W}(v,x)\in \E(W)[[x_{0}]]
\;\;\;\mbox{ for all but finitely many }g\in \Ga.$$
Thus
$$Y_{W}(Y(R_{g}u,x_{0})v,x)\in \E(W)[[x_{0}]]
\;\;\;\mbox{ for all but finitely many }g\in \Ga.$$
Since $Y_{W}$ is injective, the second assertion follows immediately.
\end{proof}

\section{Relation between $V$-modules and $V^{\otimes N}$-quasi modules}
In this section, we prove that for any $\Z$-graded vertex algebra $V$
and any positive integer $N$, $V$-module structures on a vector space
$W$ one-to-one correspond to quasimodule structures
for $V^{\otimes N}$ viewed as a $\Ga$-vertex algebra.
In view of \cite{li-tq}, this can be considered as a version of
a theorem of Barron, Dong and Mason \cite{bdm}
in terms of quasimodules.

The following is a simple result which is useful in application:

\bl{lpregphi}
Let $V$ be a $\Ga$-vertex algebra and let $(W,Y_{W})$ be a
quasimodule for $V$ viewed as a vertex algebra.
Assume that $S$ is a generating subset of $V$ as a vertex algebra such that
$\{ Y_{W}(u,x)\;|\; u\in S\}$ is $\phi(\Ga)$-local and
\begin{eqnarray}
Y_{W}(R_{g}(u),x)=Y_{W}(u,\phi(g)x)\;\;\;\mbox{ for }g\in \Ga,\; u\in S.
\end{eqnarray}
Then $(W,Y_{W})$ is a $V$-quasimodule.
\el

\begin{proof} As $(W,Y_{W})$ is a quasimodule for $V$ viewed as a
  vertex algebra,
$\bar{V}=\{ Y_{W}(v,x)\;|\; v\in V\}$ is a closed quasilocal subspace
of $\E(W)$, containing $1_{W}$, and $Y_{W}$ is a vertex-algebra 
homomorphism {}from $V$ to $(\bar{V},1_{W},Y_{\E})$.
Set $\bar{S}=\{ Y_{W}(u,x)\;|\; u\in S\}\subset \bar{V}$.
Since $S$ generates $V$ as a vertex algebra, $\bar{S}$ generates 
$\bar{V}$ as a vertex algebra. With $\bar{S}$ being $\phi(\Ga)$-local, 
{}from Theorem \ref{tga-old},  $\bar{V}$ is $\phi(\Ga)$-local. 
Set 
$$K=\{ v\in V\;|\; Y_{W}(R_{g}v,x)=Y_{W}(v,\phi(g)x)
\;\;\;\mbox{ for }g\in \Ga\}.$$
It remains to to prove $K=V$. Clearly,
$\{{\bf 1}\}\cup S\subset K$.
Let $u,v\in K,\; g\in \Ga$, $w\in W$. There exist nonzero
polynomials $f(x_{0},x)$ and $g(x_{0},x)$ such that
\begin{eqnarray*}
& &f(x_{0},x)
Y_{W}(Y(R_{g}u,\phi(g)^{-1}x_{0})R_{g}v,x)w
=f(x_{0},x)Y_{W}(R_{g}u,\phi(g)^{-1}x_{0}+x)Y_{W}(R_{g}v,x)w,\\
& &g(x_{0},x)Y_{W}(u,x_{0}+\phi(g)x)Y_{W}(v,\phi(g)x)w
=g(x_{0},x)Y_{W}(Y(u,x_{0})v,\phi(g)x)w.
\end{eqnarray*}
Then
\begin{eqnarray*}
& &f(x_{0},x)g(x_{0},x)
Y_{W}(R_{g}Y(u,x_{0})v,x)w\\
&=&f(x_{0},x)g(x_{0},x)
Y_{W}(Y(R_{g}u,\phi(g)^{-1}x_{0})R_{g}v,x)w\\
&=&f(x_{0},x)g(x_{0},x)
Y_{W}(R_{g}u,\phi(g)^{-1}x_{0}+x)Y_{W}(R_{g}v,x)w\\
&=&f(x_{0},x)g(x_{0},x)
Y_{W}(u,x_{0}+\phi(g)x)Y_{W}(v,\phi(g)x)w\\
&=&f(x_{0},x)g(x_{0},x)
Y_{W}(Y(u,x_{0})v,\phi(g)x)w,
\end{eqnarray*}
{}from which we obtain
$$Y_{W}(R_{g}Y(u,x_{0})v,x)w=Y_{W}(Y(u,x_{0})v,\phi(g)x)w.$$
It follows that $K$ is closed. As $S$ generates $V$, we must have $K=V$.
\end{proof}

Let $V$ be a $\Z$-graded vertex algebra and let $N$ be 
a fixed positive integer.
{}From \cite{fhl}, we have a tensor product vertex algebra 
$V^{\otimes N}$, which is naturally $\Z$-graded.
Let $\sigma$ be the permutation automorphism of $V^{\otimes N}$ defined by
$$\sigma(v^{1}\otimes \cdots \otimes v^{N-1}\otimes v^{N})
=v^{2}\otimes \cdots \otimes v^{N}\otimes v^{1}$$
for $v^{1},\dots,v^{N}\in V$. Set
$$\Ga=\<\sigma\>\subset \Aut (V^{\otimes N}).$$
Let $\phi: \Ga\rightarrow \C^{\times}$ be the group homomorphism 
defined by $\phi(\sigma)=\exp (2\pi i/N)$.
Denote by $\Gamma_{N}$ the group of $n$th roots of unity:
$$\Gamma_{N}=\{ \exp (2\pi i r/N)\;|\; r=0,\dots,N-1\}\subset
\C^{\times}.$$
We identify $V$ as a vertex subalgebra of
$V^{\otimes N}$ through the map
$$v\mapsto v\otimes {\bf 1}\otimes \cdots \otimes {\bf 1}.$$
Then $\sigma^{j}(V),\; j=0,\dots,N-1$, are 
mutually commuting graded vertex subalgebras and
$V^{\otimes N}=\otimes_{j=0}^{N-1}\sigma^{j}(V)$.
For $g\in \Ga$, we define
$$R_{g}a=\phi(g)^{-L(0)}g(a)\;\;\;\mbox{ for }g\in \Ga,\; a\in V^{\otimes N}.$$
In view of Example \ref{example-Zva}, $V^{\otimes N}$ 
becomes a $\Ga$-vertex algebra.
With all these notations we have:

\bt{tmodule-quasimodule}
Let $(W,\overline{Y}_{W})$ be a $V^{\otimes N}$-quasimodule.
Denote by $Y_{W}$ the restriction map of $\overline{Y}_{W}$ to $V$.
Then $(W,Y_{W})$ is a $V$-module.
On the other hand, for any $V$-module $(W,Y_{W})$,
the linear map $Y_{W}$ from $V$ to $\Hom (W,W((x)))$ can be extended
to a linear map $\overline{Y}_{W}$ from $V^{\otimes N}$ to $\Hom (W,W((x)))$
such that $(W,\overline{Y}_{W})$ is a $V^{\otimes N}$-quasimodule.
Furthermore, such an extension is unique.
\et

\begin{proof} Let $(W,\overline{Y}_{W})$ be a $V^{\otimes N}$-quasimodule.
Then $(W,Y_{W})$ is a quasimodule for $V$ (a vertex algebra).
For $u,v\in V,\; g\in \Ga$ with $g\ne 1$, we have
$$Y(g(u),x)v\in (V^{\otimes N})[[x]].$$
Using this property and Proposition \ref{pmore1} we get
\begin{eqnarray*}
& &[\overline{Y}_{W}(u,x_{1}),\overline{Y}_{W}(v,x_{2})]\\
&=&\sum_{g\in \Ga}\Res_{x_{0}}
x_{1}^{-1}\delta\left(\frac{\phi(g) x_{2}+x_{0}}{x_{1}}\right)
\overline{Y}_{W}(Y(\phi(g)^{-L(0)}g(u),\phi(g)^{-1} x_{0})v,x_{2})\\
&=&\Res_{x_{0}}
x_{1}^{-1}\delta\left(\frac{x_{2}+x_{0}}{x_{1}}\right)
\overline{Y}_{W}(Y(u,x_{0})v,x_{2}).
\end{eqnarray*}
It follows that
$$(x_{1}-x_{2})^{k}
[\overline{Y}_{W}(u,x_{1}),\overline{Y}_{W}(v,x_{2})]=0$$
for any nonnegative integer $k$ with
$x^{k}Y(u,x)v\in V[[x]]$. In view of Lemma \ref{lsimple2},
$(W,Y_{W})$ is a $V$-module.

For the other direction, first we prove the uniqueness.
Suppose that a $V$-module structure $Y_{W}$ on $W$
is extended to a $V^{\otimes N}$-quasimodule structure 
$\overline{Y}_{W}$.
By (\cite{li-ga}, Proposition 5.4), $\overline{Y}_{W}$ 
is a vertex-algebra homomorphism.
For $g\in \Ga,\; v\in V$, we have
\begin{eqnarray}
\ \ \ \ \overline{Y}_{W}(gv,x)
=\overline{Y}_{W}(\phi(g)^{L(0)}v,\phi(g)x)
=Y_{W}(\phi(g)^{L(0)}v,\phi(g)x).
\end{eqnarray}
Since $g(V)$, $g\in \Ga$, generate $V^{\otimes N}$ as a vertex algebra,
the uniqueness is clear.

Now it remains to establish the existence. For $g\in \Ga$, set
$$U[g]=\{ Y_{W}(\phi(g)^{L(0)}v,\phi(g) x)\;|\; v\in V\}
\subset \Hom (W,W((x)))=\E(W).$$
Furthermore, set
$$U=\sum_{g\in \Ga}U[g] \subset \E(W).$$
For $u,v\in V,\; \alpha,\beta\in \C^{\times}$, 
there exists a nonnegative integer $k$ such that
$$(x_{1}-x_{2})^{k}
[Y_{W}(\alpha^{L(0)}u,x_{1}),Y_{W}(\beta^{L(0)}v,x_{2})]=0.$$
Consequently,
\begin{eqnarray}\label{egamma-local}
(x_{1}-\alpha^{-1}\beta x_{2})^{k}
[Y_{W}(\alpha^{L(0)}u,\alpha x_{1}),Y_{W}(\beta^{L(0)}v,\beta
  x_{2})]=0.
\end{eqnarray}
It follows that $U$ is $\Gamma_{N}$-local.
By Theorem \ref{tga-old}, $U$ generates a 
$\Gamma_{N}$-vertex algebra $\<U\>_{\Gamma_{N}}$ with 
$W$ as a $\Gamma_{N}$-quasimodule.
Furthermore, for every $g\in \Ga$, {}from (\ref{egamma-local})
$U[g]$ is local, so
{}from \cite{li-local}, $U[g]$ generates
a vertex algebra with $W$ as a module, where for $u,v\in V$,
\begin{eqnarray*}
& &Y_{\E}(Y_{W}(\phi(g)^{L(0)}u,\phi(g) x),x_{0})
Y_{W}(\phi(g)^{L(0)}v,\phi(g) x)\\
&=&\Res_{x_{1}}x_{0}^{-1}\delta\left(\frac{x_{1}-x}{x_{0}}\right)
Y_{W}(\phi(g)^{L(0)}u,\phi(g) x_{1})Y_{W}(\phi(g)^{L(0)}v,\phi(g) x)\\
& &-\Res_{x_{1}}x_{0}^{-1}\delta\left(\frac{x-x_{1}}{-x_{0}}\right)
Y_{W}(\phi(g)^{L(0)}v,\phi(g) x)Y_{W}(\phi(g)^{L(0)}u,\phi(g) x_{1})\\
&=&\Res_{x_{1}}x_{1}^{-1}\delta\left(\frac{x+x_{0}}{x_{1}}\right)
Y_{W}(Y(\phi(g)^{L(0)}u,\phi(g) x_{0})\phi(g)^{L(0)}v,\phi(g) x)\\
&=&Y_{W}(\phi(g)^{L(0)}Y(u,x_{0})v,\phi(g) x).
\end{eqnarray*}
This shows that the linear map 
$$f_{g}: V\rightarrow U[g]\subset
\<U\>_{\Gamma},\;\; v\mapsto Y_{W}(\phi(g)^{L(0)}v,\phi(g)x)$$
 is a vertex-algebra homomorphism.
Furthermore, if $g\ne h$, with the relation
(\ref{egamma-local}) it follows from
(\cite{li-ga}, Proposition 4.8) that 
for $a(x)\in U[g],\; b(x)\in U[h]$,
$$[Y_{\E}(a(x),x_{1}),Y_{\E}(b(x),x_{2})]=0.$$
Then vertex-algebra homomorphisms $f_{g}$ $(g\in \Ga)$
give rise to a vertex-algebra homomorphism 
 $\overline{Y}_{W}$ {}from $V^{\otimes N}$ to $\<U\>_{\Gamma_{N}}$. 
Consequently, $W$ is a quasimodule for $V^{\otimes N}$ 
viewed as a vertex algebra, where
for $v\in V,\; g\in \Ga$,
$$\overline{Y}_{W}(g(v),x)=Y_{W}(\phi(g)^{L(0)}v,\phi(g) x).$$
For $g,h\in \Ga,\; v\in V$, we have
\begin{eqnarray*}
\overline{Y}_{W}(\phi(g)^{-L(0)}g(h(v)),x)
&=&\overline{Y}_{W}(\phi(g)^{-L(0)}gh(v),x)\\
&=&Y_{W}(\phi(h)^{L(0)}v, \phi(gh) x)\\
&=&\overline{Y}_{W}(h(v),\phi(g) x).
\end{eqnarray*}
By Lemma \ref{lpregphi}, $W$ is a quasimodule for $V^{\otimes N}$
viewed as a $\Ga$-vertex algebra.
\end{proof}

\section{Certain generalizations of twisted affine Lie algebras}
In this section we extend the results of 
Golenishcheva-Kutuzova and Kac (see \cite{gkk}) on
a certain generalization of the construction of twisted affine Lie algebras.
We show that restricted modules for such generalized twisted affine Lie algebras
are quasimodules for vertex algebras associated with (untwisted) affine Lie algebras.
We also formulate a notion of 
$\Gamma$-conformal Lie algebra, which extends the notion of 
$\Gamma$-conformal Lie algebra of \cite{gkk}. 

First we prove the following simple result:

\bl{lgkk}
Let $K$ be a Lie algebra equipped with a (possibly zero)
symmetric invariant bilinear form $\<\cdot,\cdot\>$.
Assume that a group $\Ga$ acts on $K$ 
by automorphisms preserving the bilinear form such that for any $u,v\in K$, 
$$[gu,v]=0,\ \ \ \ \<gu,v\>=0\ \ \ 
\mbox{ for all but finitely many }g\in \Ga.$$
Define a new multiplicative operation $[\cdot,\cdot]_{\Ga}$ on $K$ by
\begin{eqnarray}
[u,v]_{\Ga}=\sum_{g\in \Ga}[gu,v]
\end{eqnarray}
for $u,v\in K$. Then the subspace linearly spanned by the vectors
$gu-u$ for $g\in \Ga,\; u\in K$ is a two-sided ideal of the new
nonassociative algebra and the quotient algebra which we denote by 
$K/\Gamma$ is a Lie algebra. Define a bilinear form 
$\<\cdot,\cdot\>_{\Ga}$ on $K$ by
$$\<u,v\>_{\Ga}=\sum_{g\in \Ga}\<gu,v\>\;\;\;\mbox{ for }u,v\in K.$$
Then $\<\cdot,\cdot\>_{\Ga}$ naturally gives rise to a symmetric invariant
bilinear form on $K/\Ga$.
\el

\begin{proof} Let $g\in \Ga,\; u,v,w\in K$. We have
$$[gu-u,v]_{\Ga}=\sum_{h\in\Ga}[hgu,v]-\sum_{h\in \Ga}[hu,v]
=\sum_{k\in\Ga}[ku,v]-\sum_{h\in \Ga}[hu,v]=0.$$
Using the assumption that $\Ga$ acts on $K$ by automorphisms, we have
$$[v,u]_{\Ga}=\sum_{g\in \Ga}[g^{-1}v,u]=\sum_{g\in \Ga}g^{-1}[v,gu]
=-\sum_{g\in \Ga}g^{-1}[gu,v],$$
{}from which we get 
$$[u,v]_{\Ga}+[v,u]_{\Ga}=\sum_{g\in \Ga}([gu,v]-g^{-1}[gu,v]).$$ 
Furthermore, we have
\begin{eqnarray*}
& &[u,[v,w]_{\Ga}]_{\Ga}-[v,[u,w]_{\Ga}]_{\Ga}
-[[u,v]_{\Ga},w]_{\Ga}\\
&=&\sum_{g,h\in \Ga}[gu,[hv,w]]-
\sum_{g,h\in \Ga}[hv,[gu,w]]-
\sum_{g,h\in \Ga}[[hgu,hv],w]=0.
\end{eqnarray*}
{}From these the first assertion follows.

For the second assertion, for $g\in \Ga,\; u,v\in K$, we have
$$\<gu-u,v\>_{\Ga}=\sum_{h\in \Ga}\<h(gu-u),v\>=
\sum_{h\in \Ga}\<hgu,v\>-\sum_{h\in \Ga}\<hu,v\>=0.$$
The bilinear form $\<\cdot,\cdot\>_{\Ga}$ is symmetric as
$$\<v,u\>_{\Ga}=\sum_{g\in \Ga}\<gv,u\>
=\sum_{g\in \Ga}\<g^{-1}u,v\>=\<u,v\>_{\Ga}.$$
For $h\in \Ga,\; u,v,w\in K$, we have
$$\<hu,hv\>_{\Ga}=\sum_{g\in \Ga}\<ghu,hv\>=\sum_{g\in
  \Ga}\<h^{-1}ghu,v\>
=\<u,v\>_{\Ga}$$
and
$$\<[u,v]_{\Ga},w\>_{\Gamma}=\sum_{g,h\in \Ga}\<g[hu,v],w\>
=\sum_{g,h\in \Ga}\<[ghu,gv],w\>=\sum_{g,h\in \Ga}\<ghu,[gv,w]\>
=\<u,[v,w]_{\Ga}\>_{\Ga}.$$
{}From these the second assertion follows.
\end{proof}

\br{rinv-coinv}
{\em Let $\Ga$ be a finite group acting on a Lie algebra $K$ by automorphisms.
Then the assumption in Lemma \ref{lgkk} automatically holds.
On the one hand, we have the $\Ga$-invariant Lie subalgebra $K^{\Ga}$ 
(the set of $\Ga$-fixed points) and on the other hand,
we have the Lie algebra $K/\Ga$.
It is straightforward to show that
the linear map $\psi: K\rightarrow K^{\Ga}$ defined by
$\psi(u)=\sum_{g\in \Ga}gu$ gives rise to a Lie algebra isomorphism 
{}from $K/\Ga$ onto $K^{\Ga}$. }
\er

\br{rgamma-cl}
{\em In \cite{gkk}, Golenishcheva-Kutuzova and Kac studied 
a notion of $\Gamma$-conformal algebra and they proved that
a $\Gamma$-conformal algebra structure on a vector space $\g$ 
exactly amounts to a Lie algebra structure 
together with a group action of $\Ga$ on $\g$
by automorphisms such that for any $u,v\in \g$, 
$[gu,v]=0$ for all but finitely many $g\in \Gamma$. Furthermore,
a loop-like (or current-like) Lie algebra was associated to every
$\Gamma$-conformal algebra together with a group homomorphism 
{}from $\Ga$ to $\C^{\times}$.}
\er

The following proposition extends a result of \cite{gkk} 
with a central extension included (with a different proof):

\bp{paffine}
Let $\g$ be a (possibly infinite-dimensional) 
Lie algebra equipped with a symmetric invariant bilinear 
form $\<\cdot,\cdot\>$,
let $\Gamma$ be a subgroup of $\Aut (\g,\<\cdot,\cdot\>)$
and let $\phi$ be any group homomorphism from $\Gamma$ to $\C^{\times}$.
Assume that for $a,b\in \g$, 
\begin{eqnarray}
[ga,b]=0\;\;\mbox{ and }\; \<ga,b\>=0\;\;\;
\mbox{ for all but finitely many }g\in \Gamma. 
\end{eqnarray}
Define a bilinear multiplicative operation $[\cdot,\cdot]_{\Ga}$ on
the vector space $\g\otimes \C[t,t^{-1}]\oplus \C {\bf k}$ by
\begin{eqnarray}
[a\otimes t^{m}+\alpha {\bf k},b\otimes t^{n}+\beta {\bf k}]_{\Ga}
=\sum_{g\in \Ga}\phi(g)^{m}\left( [ga,b]\otimes t^{m+n}
+m\<ga,b\>\delta_{m+n,0}{\bf k}\right)
\end{eqnarray}
for $a,b\in \g,\; m,n\in \Z,\; \alpha,\beta\in \C$. 
Then the subspace linearly spanned 
by the elements
\begin{eqnarray}
\phi(g)^{m}(ga\otimes t^{m})-(a\otimes t^{m})
\;\;\;\mbox{ for }g\in \Ga,\; a\in \g,\; m\in \Z
\end{eqnarray}
is a two-sided ideal of the nonassociative algebra and the quotient algebra
which we denote by $\hat{\g}[\Gamma]$ is a Lie algebra.
\ep

\begin{proof} Associated to the pair $(\g,\<\cdot,\cdot\>)$,
we have the (untwisted) affine Lie algebra 
\begin{eqnarray}
\hat{\g}=\g\otimes \C[t,t^{-1}]\oplus \C {\bf k},
\end{eqnarray}
where ${\bf k}$ is central, and for $a,b\in \g,\; m,n\in \Z$,
\begin{eqnarray}
[a\otimes t^{m},b\otimes t^{n}]
=[a,b]\otimes t^{m+n}+m\delta_{m+n,0}\<a,b\>{\bf k}.
\end{eqnarray}
Let $\Ga$ act on $\hat{\g}$ by
$$g(a\otimes t^{m}+\beta {\bf k})=\phi(g)^{m}(ga\otimes t^{m})+\beta{\bf k}$$
for $g\in \Ga,\; a\in \g,\; m\in \Z,\; \beta\in \C$.
It is straightforward to see that $\Ga$ acts on 
$\hat{\g}$ by automorphisms. Furthermore,
for any $a,b\in \g,\; m,n\in \Z,\; \alpha,\beta\in \C$, we have
$$[g(a\otimes t^{m}+\alpha {\bf k}), b\otimes t^{n}+\beta {\bf k}]
=\phi(g)^{m}\left([ga,b]\otimes t^{m+n}
+m\delta_{m+n,0}\<ga,b\>{\bf k}\right)=0$$
for all but finitely many $g\in \Ga$. 
Now it follows immediately from Lemma \ref{lgkk}
with $K=\hat{\g}$.
\end{proof}

\br{rtwisted-affine}
{\em Let $\g, \<\cdot,\cdot\>$ be given as in Proposition \ref{paffine}, and let 
$\sigma$ be an order $N$ automorphism of $\g$, preserving the bilinear form
$\<\cdot,\cdot\>$. Extend $\sigma$ to an automorphism of Lie algebra $\hat{\g}$ by
$$\sigma (u\otimes t^{m}+\alpha {\bf k})=\exp (-2m\pi i/N)(\sigma (u)\otimes t^{m})
+\alpha {\bf k}$$
for $u\in \g,\; m\in \Z,\; \alpha\in \C$.
The twisted affine Lie algebra $\hat{\g}[\sigma]$ (see \cite{kac1}) 
can be realized as
the $\sigma$-fixed point subalgebra of $\hat{\g}$.
Set $\Ga=\<\sigma\>\subset \Aut(\g,\<\cdot,\cdot\>)$
and let $\phi$ be the group embedding of $\Ga$ into $\C^{\times}$ defined by
$\phi(\sigma^{n})=\exp (-2n\pi i/N)$. We let $\Ga$ act on $\hat{\g}$ as 
in the proof of Proposition \ref{paffine}. Clearly, $\hat{\g}[\sigma]$ 
is also the $\Ga$-invariant subalgebra.
In view of this and Remark \ref{rinv-coinv}, Lie algebras $\hat{\g}[\Ga]$ 
are generalizations of twisted affine Lie algebras. }
\er

\br{rco-affine}
{\em  Let $\g, \<\cdot,\cdot\>, \Ga,\; \phi$ 
be given as in Proposition \ref{paffine}.
Set $H=\ker \phi\subset \Ga$, a normal subgroup. 
In view of Lemma \ref{lgkk}, we have a 
Lie algebra $\g/H$ equipped with a symmetric invariant bilinear form
$\<\cdot,\cdot\>_{H}$. Then we have the (untwisted) 
affine Lie algebra $\widehat{\g/H}$. On the other hand,
$\Ga/H$ naturally acts on the Lie algebra $\g/H$ by automorphisms
and $\phi$ reduces to a group embedding of $\Ga/H$ into $\C^{\times}$.
In view of Proposition \ref{paffine}, we have a Lie algebra $\widehat{(\g/H)}[\Ga/H]$. }
\er

We have:

\bp{pisomorphism}
Let $\g,\<\cdot,\cdot\>,\Ga$ and $\phi$ be given as in Proposition \ref{paffine}
and set $H=\ker \phi\subset \Ga$. The Lie algebra $\hat{\g}[\Ga]$ 
is isomorphic to the Lie algebra $\widehat{(\g/H)}[\Ga/H]$.
\ep

\begin{proof} It is straightforward.
\end{proof}

\br{raffine-voa}
{\em We here review the $\Z$-graded vertex algebras
associated (untwisted) affine Lie algebras. Let $\g$ be a Lie algebra 
equipped with a symmetric invariant bilinear form $\<\cdot,\cdot\>$
and let $\hat{\g}$ be the associated affine Lie algebra.
For $a\in \g$, form the generating function
$$a(x)=\sum_{n\in \Z}a(n)x^{-n-1},$$
where $a(n)$ is an alternative notation for $a\otimes t^{n}$.
Lie algebra $\hat{\g}$ is naturally $\Z$-graded 
$\hat{\g}=\coprod_{n\in \Z}\hat{\g}_{(n)}$,
where 
$$\hat{\g}_{(0)}=\g\oplus \C {\bf k},\;\;\; \hat{\g}_{(n)}=\g\otimes t^{-n}
\;\;\;\mbox{ for }n\ne 0.$$
Set 
$$\hat{\g}_{(\le 0)}=\g\otimes \C[t]\oplus \C {\bf k},\;\;\;
\hat{\g}_{(+)}=\g\otimes t^{-1}\C[t^{-1}].$$
Let $\ell$ be a complex number and let $\C_{\ell}=\C$ be the $1$-dimensional 
$\hat{\g}_{(\le 0)}$-module with $\g\otimes \C[t]$ acting trivially and with
${\bf k}$ acting as scalar $\ell$. Form the induced module
$$V_{\hat{\g}}(\ell,0)=U(\hat{\g})\otimes_{U(\hat{\g}_{(\le 0)})} \C_{\ell},$$
which is naturally an $\N$-graded $\hat{\g}$-module (of level $\ell$).
Set ${\bf 1}=1\otimes 1$ and identify $\g$ as a subspace through the map
$a\mapsto a(-1){\bf 1}$. In fact, $\g$ is exactly the degree-one subspace.
It was known (cf. [FZ], [Lia], [LL]) that there exists a (unique) 
vertex algebra structure on $V_{\hat{\g}}(\ell,0)$
with ${\bf 1}$ as the vacuum vector and with $Y(a,x)=a(x)$ for $a\in \g$.
Furthermore, vertex algebra $V_{\hat{\g}}(\ell,0)$
satisfies the following universal property (cf. \cite{pr}):
Let $V$ be a vertex algebra and let $f$ be a linear map {}from $\g$ to
$V$ such that for $a,b\in \g$,
$$f(a)_{0}f(b)=f([a,b]),\ \ \ \ f(a)_{1}f(b)=\ell \<a,b\> {\bf 1},
\;\;\mbox{ and }\; f(a)_{n}f(b)=0\;\;\;\mbox{ for }n\ge 2.$$
Then $f$ can be extended (uniquely) to a vertex-algebra homomorphism {}from
$V_{\hat{\g}}(\ell,0)$ to $V$.
Let $\sigma\in \Aut(\g,\<\cdot,\cdot\>)$, that is,
$\sigma$ is an automorphism of Lie algebra $\g$, which preserves the
bilinear form. Then $\sigma$ extends canonically to
an automorphism of the $\Z$-graded vertex algebra
$V_{\hat{\g}}(\ell,0)$. For any subgroup $\Ga$ of 
$\Aut(\g,\<\cdot,\cdot\>)$ and any group homomorphism 
$\phi: \Ga\rightarrow \C^{\times}$, $V_{\hat{\g}}(\ell,0)$ 
is a $\Ga$-vertex algebra with $R_{g}=\phi(g)^{-L(0)}g$ for $g\in \Ga$.}
\er

The following is the main result of this section:

\bt{tmain2}
Let $\g$ be a Lie algebra equipped with a symmetric invariant bilinear 
form $\<\cdot,\cdot\>$. Let $\Gamma$ be a subgroup of 
$\Aut(\g,\<\cdot,\cdot\>)$ such that for $a,b\in \g$,
$$[ga,b]=0\;\;\;\mbox{ and }\;\;\; \<ga,b\>=0\;\;\;
\mbox{ for all but finitely many } g\in \Gamma.$$
Let $\phi$ be any group homomorphism from $\Gamma$ to $\C^{\times}$
and set $H=\ker \phi\subset \Gamma$. 
Then any restricted module $W$ of level $\ell$ 
for Lie algebra $\hat{\g}[\Ga]$ is a quasimodule
for $V_{\widehat{\g/H}}(\ell,0)$ viewed as a $\Gamma$-vertex algebra
with $Y_{W}(\bar{a},x)=a_{W}(x)$ for $a\in \g$, where $\bar{a}$ denotes 
the image of $a$ in $\g/H$ under the natural quotient map.
On the other hand, 
any quasimodule $(W,Y_{W})$ for $V_{\widehat{(\g/H)}}(\ell,0)$ viewed as a
$\Gamma$-vertex algebra is a
restricted module of level $\ell$ for Lie algebra
$\hat{\g}[\Gamma]$ with $a_{W}(x)=Y_{W}(\bar{a},x)$ for $a\in \g$.
\et

\begin{proof} Let $W$ be a restricted $\hat{\g}[\Ga]$-module 
of level $\ell$. For $a\in \g$, form the generating function
$$a_{W}(x)=\sum_{n\in \Z}\overline{a(n)}x^{-n-1}\in \E(W),$$
where $\overline{a(n)}$ denotes the operator on $W$,
associated to the image of $a\otimes t^{n}$ in $\hat{\g}[\Gamma]$.
Set
$$U=\{ a_{W}(x)\;|\; a\in \g\}\subset \E(W).$$
For $a,b\in \g$, we have
\begin{eqnarray}
[a_{W}(x_{1}),b_{W}(x_{2})]
=\sum_{g\in \Ga}[ga,b]_{W}(x_{2})
x_{1}^{-1}\delta\left(\frac{\phi(g)x_{2}}{x_{1}}\right)
+\ell \<ga,b\>\frac{\partial}{\partial x_{2}}
x_{1}^{-1}\delta\left(\frac{\phi(g)x_{2}}{x_{1}}\right).
\end{eqnarray}
It follows that $U$ is $\phi(\Ga)$-local. By Theorem \ref{tga-old}
$U$ generates a $\phi(\Ga)$-vertex algebra $\<U\>_{\phi(\Ga)}$
with $W$ as a faithful quasimodule where
$Y_{W}(\alpha(x),x_{0})=\alpha(x_{0})$ for $\alpha(x)\in \<U\>_{\phi(\Ga)}$.
As for $h\in H=\ker \phi$, 
$$(ha)_{W}(x)=\phi(h)a_{W}(\phi(h)x)=a_{W}(x),$$
we have a linear map from $\g/H$ into $\<U\>_{\phi(\Ga)}$, sending
$\bar{a}$ to $a_{W}(x)$ for $a\in \g$. Furthermore, for $a,b\in \g$ we have
\begin{eqnarray*}
& &[Y_{W}(a_{W}(x),x_{1}),Y_{W}(b_{W}(x),x_{2})]\\
&=&[a_{W}(x_{1}),b_{W}(x_{2})]\\
&=&\sum_{g\in \Ga}[ga,b]_{W}(x_{2})
x_{1}^{-1}\delta\left(\frac{\phi(g)x_{2}}{x_{1}}\right)
+\ell \<ga,b\>\frac{\partial}{\partial x_{2}}
x_{1}^{-1}\delta\left(\frac{\phi(g)x_{2}}{x_{1}}\right)\\
&=&\sum_{g\in \Ga}Y_{W}([ga,b],x_{2})
x_{1}^{-1}\delta\left(\frac{\phi(g)x_{2}}{x_{1}}\right)
+\ell \<ga,b\>\frac{\partial}{\partial x_{2}}
x_{1}^{-1}\delta\left(\frac{\phi(g)x_{2}}{x_{1}}\right).
\end{eqnarray*}
Comparing this with Proposition \ref{pmore1}, we obtain
\begin{eqnarray*}
& &a_{W}(x)_{0}b_{W}(x)=\sum_{h\in H} [ha,b]_{W}(x),\ \ \ 
a_{W}(x)_{1}b_{W}(x)=\ell \sum_{h\in H}\<ha,b\>1_{W},\\
& &a_{W}(x)_{n}b_{W}(x)=0\;\;\;\mbox{ for }n\ge 2.
\end{eqnarray*}
In view of Remark \ref{raffine-voa}, there exists a (unique) vertex-algebra
homomorphism {}from $V_{\widehat{\g/H}}(\ell,0)$ to $\<U\>_{\phi(\Ga)}$,
sending $\bar{a}$ to $a_{W}(x)$ for $a\in \g$. 
It follows from Lemma \ref{lpregphi} that $W$ is a quasimodule for
$V_{\widehat{(\g/H)}}(\ell,0)$ viewed as a vertex algebra
with $Y_{W}(\bar{a},x)=a_{W}(x)$ for $a\in \g$.
For $g\in \Ga,\; a\in \g$, we have
$$Y_{W}(R_{g}a,x)=\phi(g)^{-1}Y_{W}(ga,x)=\phi(g)^{-1}(ga)_{W}(x)
=a_{W}(\phi(g)x)=Y_{W}(a,\phi(g)x).$$
Now it follows from Lemma \ref{lpregphi} that $W$ is a quasimodule for
$V_{\widehat{(\g/H)}}(\ell,0)$ viewed as a $\Ga$-vertex algebra.

On the other hand, let $(W,Y_{W})$ be a quasimodule for
$V_{\widehat{(\g/H)}}(\ell,0)$ viewed as a $\Gamma/H$-vertex algebra. 
For $a\in \g$, set $a_{W}(x)=Y_{W}(\bar{a},x)$.
For $g\in \Ga,\; a\in \g$, we have
$$(ga)_{W}(x)=Y_{W}(\overline{ga},x)=Y_{W}(\bar{g}\bar{a},x)
=\bar{\phi}(\bar{g})Y_{W}(\bar{a},\bar{\phi}(\bar{g})x)
=\phi(g)a_{W}(\phi(g)x).$$
Notice that for $u,v\in \g$, with $\bar{u},\bar{v}\in V_{\widehat{(\g/H)}}(\ell,0)$, we have
$$\bar{u}_{0}\bar{v}=[\bar{u},\bar{v}]=\sum_{h\in H}\overline{[hu,v]},
\;\; \bar{u}_{1}\bar{v}=\ell \<\bar{u},\bar{v}\>{\bf 1}=\sum_{h\in H}\ell \<hu,v\>{\bf 1},
\;\; \bar{u}_{n}\bar{v}=0
\;\;\;\mbox{ for }n\ge 2.$$
For $a,b\in \g$, using Proposition \ref{pmore1} we have
\begin{eqnarray*}
& &[a_{W}(x_{1}),b_{W}(x_{2})]\\
&=&[Y_{W}(\bar{a},x_{1}),Y_{W}(\bar{b},x_{2})]\\
&=&\sum_{\bar{g}\in \Ga/H}\Res_{x_{0}}
x_{1}^{-1}\delta\left(\frac{\bar{\phi}(\bar{g}) x_{2}+x_{0}}{x_{1}}\right)
Y_{W}(Y(R_{\bar{g}}\bar{a},\bar{\phi}(\bar{g})^{-1}x_{0})\bar{b},x_{2})\\
&=& \sum_{\bar{g}\in \Ga/H}Y_{W}([\bar{g}\bar{a},\bar{b}],x_{2})
x_{1}^{-1}\delta\left(\frac{\bar{\phi}(\bar{g}) x_{2}}{x_{1}}\right)
+\ell \<\bar{g}\bar{a},\bar{b}\>
\frac{\partial}{\partial x_{2}}
x_{1}^{-1}\delta\left(\frac{\bar{\phi}(\bar{g}) x_{2}}{x_{1}}\right)\\
&=& \sum_{g\in \Ga}Y_{W}(\overline{[ga,b]},x_{2})
x_{1}^{-1}\delta\left(\frac{\phi(g) x_{2}}{x_{1}}\right)+\ell \<ga,b\>
\frac{\partial}{\partial x_{2}}
x_{1}^{-1}\delta\left(\frac{\phi(g) x_{2}}{x_{1}}\right)\\
&=& \sum_{g\in \Ga}[ga,b]_{W}(x_{2})
x_{1}^{-1}\delta\left(\frac{\phi(g) x_{2}}{x_{1}}\right)+\ell \<ga,b\>
\frac{\partial}{\partial x_{2}}
x_{1}^{-1}\delta\left(\frac{\phi(g) x_{2}}{x_{1}}\right).
\end{eqnarray*}
It follows that $W$ is a restricted $\hat{\g}[\Ga]$-module of level $\ell$
with $u_{W}(x)=Y_{W}(\bar{u},x)$ for $u\in \g$.
\end{proof}

\bex{general-example}
{\em Let $\Ga$ be a group as before.
We define an associative algebra $gl_{\Ga}$ with 
a basis $\{ E_{\alpha,\beta}\;|\; \alpha,\beta\in \Ga\}$ such that
\begin{eqnarray}
E_{\alpha,\beta}\cdot E_{\lambda,\mu}=\delta_{\beta,\lambda}E_{\alpha,\mu}
\;\;\;\mbox{ for }\alpha,\beta,\lambda,\mu\in \Ga.
\end{eqnarray}
Equip $gl_{\Ga}$ with a bilinear form $\<\cdot,\cdot\>$ defined by
\begin{eqnarray}
\<E_{\alpha,\beta},E_{\lambda,\mu}\>=\delta_{\alpha,\mu}\delta_{\beta,\lambda}
\;\;\;\mbox{ for }\alpha,\beta,\lambda,\mu\in \Ga.
\end{eqnarray}
Clearly, this form is nondegenerate, symmetric and associative (invariant).
For any associative algebra $A$ equipped with a nondegenerate symmetric 
invariant bilinear form $\<\cdot,\cdot\>$, the tensor product associative algebra
$A\otimes gl_{\Ga}$ has a nondegenerate symmetric invariant bilinear form with
\begin{eqnarray}
\< a\otimes E_{\alpha,\beta},b\otimes E_{\lambda,\mu}\>
=\<a,b\>\<E_{\alpha,\beta},E_{\lambda,\mu}\>
=\delta_{\alpha,\mu}\delta_{\beta,\lambda}\<a,b\>
\end{eqnarray}
for $a,b\in A,\; \alpha,\beta,\lambda,\mu\in \Ga$. The bilinear form 
is still invariant with $A\otimes gl_{\Ga}$ viewed as a Lie algebra.
Then we have a (generalized) affine Lie algebra
\begin{eqnarray}
\widehat{A\otimes gl_{\Ga}}
=(A\otimes gl_{\Ga})\otimes \C[t,t^{-1}]\oplus \C {\bf k}.
\end{eqnarray}
Let $\Gamma$ act on $A\otimes gl_{\Gamma}$ by
\begin{eqnarray}
T_{g}(a\otimes E_{\alpha,\beta})=a\otimes E_{g\alpha,g\beta}
\;\;\;\mbox{ for }g,\alpha,\beta\in \Gamma,\; a\in A.
\end{eqnarray}
Clearly, this defines an action of $\Gamma$ on $A\otimes gl_{\Gamma}$ 
by automorphisms
and for $g\in \Gamma$, $T_{g}$ preserves the bilinear form.
Let $\phi$ be any group homomorphism from $\Gamma$ to $\C^{\times}$.
Set 
\begin{eqnarray}
R_{g}=\phi(g)^{-1}T_{g}\;\;\;\mbox{ for }g\in \Gamma.
\end{eqnarray}
It is clear that for any $u,v\in A\otimes gl_{\Gamma}$,
\begin{eqnarray}
[R_{g}u,v]=0\;\;\;\mbox{and  }\;\<R_{g}u,v\>=0
\;\;\;\mbox{ for all but finitely many }g\in \Gamma.
\end{eqnarray}
In view of Proposition \ref{paffine} we have a Lie algebra
$\widehat{A\otimes gl_{\Gamma}}[\Gamma]$. By Theorem \ref{tmain2}, 
any restricted module $W$ of level $\ell$ for 
$\widehat{A\otimes gl_{\Gamma}}[\Gamma]$
is naturally a quasimodule for some vertex algebra.
Now, let $\Gamma=\Z^{k}$ be a free abelian group of rank $k$ and let
$\hbar=(h_{1},\dots,h_{k})\in \R^{k}$. Define a group homomorphism
$\phi_{\hbar}$ from $\Z^{k}$ to $\C^{\times}$ by
\begin{eqnarray}
\phi_{\hbar}(n_{1},\dots,n_{k})
=e^{\pi i(h_{1}n_{1}+\cdots +h_{k}n_{k})}.
\end{eqnarray}
We have a Lie algebra $\widehat{gl_{\Z^{k}}}[\Z^{k}]$.
This is the Lie algebra $\hat{A}_{\hbar}$ (with central extension)
in [GKL1.2] (cf. \cite{gkk}).}
\eex

\br{rvoa-quasimodule}
{\em Let $V$ be a vertex operator algebra in the sense of
\cite{flm} and \cite{fhl}, let $\Ga$ be a group of automorphisms of
$V$ and let $\phi: \Ga\rightarrow \C^{\times}$ be a group
homomorphism. View $V$ as a $\Gamma$-vertex algebra as in Example \ref{example-Zva}.
Let $(W,Y_{W})$ be a $V$-quasimodule.
For $g\in \Ga$, as $g(\omega)=\omega$,  we have
$$R_{g}\omega=\phi(g)^{-L(0)}g(\omega)=\phi(g)^{-2}\omega,$$
recalling that $\omega$ is the conformal vector of $V$.
Then
\begin{eqnarray*}
Y_{W}(\omega,x)=\phi(g)^{2}Y_{W}(\omega,\phi(g)x).
\end{eqnarray*}
If $\phi(g)$ is not a root of unity for some $g\in \Ga$, 
then $Y_{W}(\omega,x)=0$. Thus if $V$ is simple with $\omega\ne 0$
and if $\phi(g)$ is not a root of unity for some $g\in \Ga$,
$V$ does not have a nonzero quasimodule 
for $V$ viewed as a $\Gamma$-vertex algebra.
Nevertheless, as we show by examples in the next Remark, there are
nontrivial quasimodules for $V$ viewed just as a vertex algebra.}
\er

\br{rsln+1}
{\em  Consider the simple Lie algebra $sl_{n+1}$
as a subalgebra of $gl_{\infty}\;(=gl_{\Z})$ in the obvious way. 
Let $K$ be the linear span of $T_{m}sl_{n+1}$ for $m\in \Z$, that is,
$K$ is linearly spanned by
$$E_{i+m,j+m},\ \  E_{i+m,i+m}-E_{j+m,j+m}$$
for $1\le i\ne j\le n+1,\; m\in \Z$.
Then $K$ is a Lie subalgebra of $gl_{\infty}$ with $\Z$ as
a group of automorphisms. Let $\phi$ be any injective group homomorphism 
{}from $\Z$ to $\C^{\times}$. We have a Lie algebra
$\hat{K}[\Z]$. By Theorem \ref{tmain2}, any restricted $\hat{K}[\Z]$-module
of level $\ell\in \C$ is naturally a quasimodule for $V_{\hat{K}}(\ell,0)$.
With $sl_{n+1}$ being a Lie subalgebra of $K$, $V_{\widehat{sl_{n+1}}}(\ell,0)$ 
is naturally a vertex subalgebra of  $V_{\hat{K}}(\ell,0)$.
Consequently, any restricted $\hat{K}[\Z]$-module
of level $\ell\in \C$ is naturally a quasimodule 
for $V_{\widehat{sl_{n+1}}}(\ell,0)$
(as a vertex algebra).
More generally, for any finite-dimensional simple Lie algebra $\g$,
one can embed $\g$ into $gl_{\infty}$,
we can obtain nontrivial $V_{\hat{\g}}(\ell,0)$-quasimodules.}
\er

Next, we extend the notion of $\Ga$-conformal Lie algebra of \cite{gkk}.
First, recall that a {\em conformal Lie algebra} \cite{kac2}, 
also known as a {\em vertex Lie algebra} \cite{pr} (cf. \cite{dlm-vl}),
is a vector space $C$ equipped with a
linear operator $T$ and a linear map
\begin{eqnarray}
Y_{-}: & &C\rightarrow \Hom (C,x^{-1}C[x^{-1}])\nonumber\\
& & a\mapsto Y_{-}(a,x)=\sum_{n\ge 0}a_{n}x^{-n-1}
\end{eqnarray}
such that the following conditions hold for $a,b\in C$:
\begin{eqnarray}
& &[T,Y_{-}(a,x)]=\frac{d}{dx}Y_{-}(a,x),\\
& &Y_{-}(a,x)b={\rm Sing} \left(e^{xT}Y_{-}(b,-x)a\right),\\
& &[Y_{-}(a,x_{1}),Y_{-}(b,x_{2})]
={\rm Sing}\left(Y_{-}(Y_{-}(a,x_{1}-x_{2})b,x_{2})\right),
\end{eqnarray}
where ${\rm Sing}$ stands for the singular part.

Associated to a conformal Lie algebra $C$ 
one has a Lie algebra $\Lie (C)$ (see \cite{pr}), where 
\begin{eqnarray}
\Lie (C)=(C\otimes \C[t,t^{-1}])/
(T\otimes 1+1\otimes d/dt)(C\otimes \C[t,t^{-1}]),
\end{eqnarray}
as a vector space, and for $a,b\in C,\; m,n\in \Z$,
\begin{eqnarray}
[a\otimes t^{m},b\otimes t^{n}]
=\sum_{i\ge 0}\binom{m}{i}(a_{i}b)\otimes t^{m+n-i}.
\end{eqnarray}
Denote by $\rho$ the natural quotient map from 
$C\otimes \C[t,t^{-1}]$ onto $\Lie(C)$. For $u\in C,\; n\in \Z$, set
$$u(n)=\rho(u\otimes t^{n})\in \Lie(C)$$
and form the generating function
$$u(x)=\sum_{n\in \Z}u(n)x^{-n-1}\in \Lie(C)[[x,x^{-1}]].$$
Set
$$\Lie(C)^{+}=\rho(C\otimes \C[t]),\ \ \ \
\Lie(C)^{-}=\rho(C\otimes t^{-1}\C[t^{-1}]).$$
Then $\Lie(C)^{\pm}$ are Lie subalgebras and we have
$$\Lie(C)=\Lie(C)^{+}\oplus \Lie (C)^{-}.$$
Letting $\Lie(C)^{+}$ act trivially on $\C$, we form the induced module
$$V_{C}=U(\Lie(C))\otimes_{U(\Lie(C)^{+})}\C.$$
Set ${\bf 1}=1\otimes 1\in V_{C}$.
Identify $C$ as a subspace of $V_{C}$ through the linear map
$u\mapsto u(-1){\bf 1}$. 
There exists a unique vertex algebra structure on $V_{C}$ 
with ${\bf 1}$ as the vacuum vector and with $Y(u,x)=u(x)$ for $u\in
C$ (see \cite{pr}).
Furthermore, $C$ generates $V_{C}$ as a vertex algebra.
It was proved in \cite{pr} that for any linear map $f$ 
{}from $C$ into a vertex algebra $V$ such that
$$fT(u)=\D f(u),\;\;\;f(u_{n}v)=f(u)_{n}f(v)\;\;\;\mbox{ for }u,v\in C,\; n\ge 0,$$
$f$ can be extended uniquely to a vertex-algebra homomorphism from $V_{C}$ to $V$.

An {\em automorphism} of a conformal Lie algebra $C$ is 
a linear automorphism $\theta$ of $C$ such that
$T\theta=\theta T$ and
$\theta Y_{-}(u,x)v= Y_{-}(\theta (u),x)\theta (v)$ for $u,v\in C$.
We have the following straightforward analogue of Lemma \ref{lgkk}:

\bl{lC-modulo-H}
Let $C$ be a conformal Lie algebra and let $H$ be a group
acting on $C$ by automorphisms such that for any $a,b\in C$,
$Y_{-}(hu,x)v=0$ for all but finitely many $h\in H$.
Then the linear map 
$Y_{-}^{H}: C\rightarrow \Hom(C,x^{-1}C[x^{-1}])$, defined by
$$Y_{-}^{H}(u,x)v=\sum_{h\in H}Y_{-}(R_{h}u,x)v$$
for $u,v\in C$, naturally gives rise to a conformal Lie algebra structure
on the quotient space, denoted by  $C/H$, of $C$ modulo the subspace
 linearly spanned by the vectors
$R_{h}a-a$ for $h\in H,\; a\in C$.
\el

The following notion, which is parallel to the notion of
$\Gamma$-vertex algebra, extends the notion of 
$\Gamma$-conformal algebra in \cite{gkk}:

\bd{dgamma-cla}
{\em Let $\Gamma$ be a group as before. 
A {\em $\Gamma$-conformal Lie algebra} is a 
conformal Lie algebra $(C, Y_{-},T)$ equipped with 
group homomorphisms
\begin{eqnarray*}
& & R: \Gamma\rightarrow \GL(C);\; g\mapsto R_{g},\\
& &\phi: \Gamma\rightarrow \C^{\times}
\end{eqnarray*}
such that for any $a,b\in C$,
\begin{eqnarray}
& &TR_{g}=\phi(g)R_{g}T,\label{eTRcomm}\\
& &R_{g}Y_{-}(a,x)R_{g^{-1}}=Y_{-}(R_{g}a,\phi(g)^{-1}x),
\label{eRconjugate}\\
& &Y_{-}(R_{g}a,x)b=0\;\;\;\mbox{ for all but finitely many
}g\in \Gamma.
\end{eqnarray}}
\ed
 
Notice that in terms of components, (\ref{eRconjugate}) amounts to
$$R_{g}(u_{m}v)=\phi(g)^{m+1}(R_{g}u)_{m}R_{g}v\;\;\;\mbox{ for }m\in \Z.$$

We have the following analogue of Proposition \ref{paffine}:

\bp{pg-k-k}
Let $\Ga$ be a group and let $C$ be a $\Gamma$-conformal Lie algebra. 
Define a bilinear multiplication
on $C\otimes \C[t,t^{-1}]$ by
\begin{eqnarray}
[u\otimes t^{m},v\otimes t^{n}]_{\Ga}
=\sum_{g\in \Gamma}\sum_{i\ge 0}\binom{m}{i}\phi(g)^{m+1}
\left((R_{g}u)_{i}v\otimes t^{m+n-i}\right)
\end{eqnarray}
for $u,v\in C,\; m,n\in \Z$.
Then the subspace linearly spanned by the elements 
$$T(u)\otimes t^{m}+mu\otimes t^{m-1},\; \;
\phi(g)^{m+1}R_{g}u\otimes t^{m}-u\otimes t^{m} $$
for $g\in \Gamma,\; u\in C,\; m\in \Z$ is a two-sided ideal of the
nonassociative algebra $C\otimes \C[t,t^{-1}]$ and
the quotient algebra is a Lie algebra, which we denote by $\hat{C}[\Gamma]$.
\ep

\begin{proof} Associated to the conformal Lie algebra $C$,
we have a Lie algebra $\Lie(C)$. 
Let $\Gamma$ act on $C\otimes \C[t,t^{-1}]$ by 
$$g(u\otimes t^{m})=\phi(g)^{m+1}(R_{g}u\otimes t^{m})
\;\;\;\mbox{ for }g\in \Gamma,\; u\in C,\; m\in \Z.$$
For $g\in \Gamma,\; u,v\in C,\; m,n\in \Z$, we have
\begin{eqnarray*}
[g(u\otimes t^{m}),g(v\otimes t^{n})]
&=&\sum_{i\in \N}\binom{m}{i}\phi(g)^{m+n+2}
(R_{g}u)_{i}(R_{g}v)\otimes t^{m+n-i}\\
&=&\sum_{i\in \N}\binom{m}{i}\phi(g)^{m+n+1-i}
R_{g}(u_{i}v)\otimes t^{m+n-i}\\
&=&g[u\otimes t^{m},v\otimes t^{n}].
\end{eqnarray*}
Furthermore, using the relation (\ref{eTRcomm}) we have
\begin{eqnarray*}
g(T\otimes 1+1\otimes d/dt)(u\otimes t^{m})
&=&g(Tu\otimes t^{m}+u\otimes mt^{m-1})\\
&=&\phi(g)^{m+1}R_{g}Tu\otimes t^{m}
+\phi(g)^{m} R_{g}u\otimes mt^{m-1}\\
&=&\phi(g)^{m}TR_{g}u\otimes t^{m}
+\phi(g)^{m} R_{g}u\otimes mt^{m-1}\\
&=&\phi(g)^{-1}(T\otimes 1+d/dt)g(u\otimes t^{m}).
\end{eqnarray*}
It follows that $\Ga$ naturally acts on the Lie algebra $\Lie(C)$ 
by automorphisms. We have
\begin{eqnarray*}
\sum_{g\in \Ga}[g(u\otimes t^{m}),v\otimes t^{n}]
&=&\sum_{g\in \Ga}\phi(g)^{m+1}[R_{g}u\otimes t^{m},v\otimes t^{n}]\\
&=&\sum_{g\in \Ga}\sum_{i\in \N}\binom{m}{i}\phi(g)^{m+1}
((R_{g}u)_{i}v\otimes t^{m+n-i})\\
&=&[u\otimes t^{m},v\otimes t^{n}]_{\Ga}.
\end{eqnarray*}
Now it follows immediately from Lemma \ref{lgkk} with $K=\Lie(C)$.
\end{proof}

\bl{llift}
Let $C$ be a $\Ga$-conformal algebra and let $V_{C}$ be the associated 
vertex algebra. Then the group action of $\Ga$ on $C$ 
can be extended uniquely to a group action of $\Ga$ 
on $V_{C}$ such that $V_{C}$ 
becomes a $\Ga$-vertex algebra.
\el

\begin{proof} As $C$ generates $V_{C}$ as a vertex algebra, 
the uniqueness is clear. In the proof of Proposition \ref{pg-k-k}
we have proved that for $g\in \Ga$, the action of $g$ on
$C\otimes \C[t,t^{-1}]$, defined by
$$g(u\otimes t^{m})=\phi(g)^{m+1}(R_{g}u\otimes t^{m})
\;\;\;\mbox{ for }u\in C,\; m\in \Z,$$
reduces to an automorphism
of the associated Lie algebra $\Lie (C)$.
Clearly, $g$ preserves the polar decomposition.
Then $g$ gives rise to a linear automorphism, denoted by $R_{g}$, of
$V_{C}$ with $g({\bf 1})={\bf 1}$ and we have
\begin{eqnarray*}
R_{g}a_{m}v=\phi(g)^{m+1}(R_{g}a)_{m}R_{g}v
\;\;\;\mbox{ for }a\in C,\,v\in V_{C},\; m\in \Z.
\end{eqnarray*}
As $C$ generates $V_{C}$ as a vertex algebra, 
it follows from induction (and the Jacobi identity of 
vertex algebra $V_{C}$) that
$$R_{g}u_{m}v=\phi(g)^{m+1}(R_{g}u)_{m}R_{g}v
\;\;\;\mbox{ for all }u,v\in V_{C}.$$
It is easy to see that this defines a group action of $\Ga$ on $V_{C}$.
Therefore, $V_{C}$ is a $\Ga$-vertex algebra.
\end{proof}

Let $C$ be a $\Ga$-conformal Lie algebra. Set $H=\ker \phi\subset \Ga$.  
Notice that for any $h\in H$,
$R_{h}$ is an automorphism of conformal Lie algebra $C$.
Thus $H$ acts on $C$ by automorphisms.
By Lemma \ref{lC-modulo-H}, we have 
a conformal Lie algebra $C/H$.
It is clear that $C/H$ with the natural $\Ga/H$-action is 
also a $\Ga/H$-conformal Lie algebra.
By Lemma \ref{llift}, $V_{C/H}$ is naturally a $\Ga/H$-vertex algebra.
We define a notion of restricted module for the Lie algebra $\hat{C}[\Ga]$ 
in the obvious way and for a restricted module $W$ 
we define the notion $a_{W}(x)$ for $a\in C$ in the obvious way.
With all these, by slightly modifying the proof of
Theorem \ref{tmain2}  we have:

\bt{tmain1}
Let $\Gamma$ be a group and let $C$ be a
$\Gamma$-conformal Lie algebra. Set $H=\ker \phi\subset \Ga$. Then
any restricted module $W$ for the Lie
algebra $\hat{C}[\Gamma]$ is naturally a $V_{C/H}$-quasimodule
with $Y_{W}(a,x)=a_{W}(x)$ for $a\in C$.
On the other hand, any $V_{C/H}$-quasimodule $W$ is naturally a
restricted module for the Lie algebra $\hat{C}[\Ga]$
with $a_{W}(x)=Y_{W}(a,x)$ for $a\in C$.
\et

\end{document}